\documentclass[sigconf, authorversion, nonacm]{acmart}

\usepackage{balance} 
\usepackage{bm}
\usepackage{amsmath}
\usepackage{multirow}

\DeclareMathOperator*{\argmin}{arg\,min}
\usepackage[T1]{fontenc}

\title{Electric Vehicle Routing Problem for Emergency Power Supply: Towards Telecom Base Station Relief}
\titlenote{Preprint. Work in progress.}

\author{Daisuke Kikuta$^{\dagger}$, Hiroki Ikeuchi$^{\dagger}$, Kengo Tajiri$^{\dagger}$, \\
Yuta Toyama$^{\ddagger}$, Masaski Nakamura$^{\ddagger}$, Yuusuke Nakano$^{\dagger}$}
\affiliation{
  \institution{${}^{\dagger}$NTT Corporation, ${}^{\ddagger}$NTT DOCOMO \\
               \texttt{\{daiuske.kikuta,hiroki.ikeuchi,kengo.tajiri,yuusuke.nakano\}@ntt.com} \\ 
               \texttt{\{yuuta.toyhama.zg,masaki.nakamura\}@nttdocomo.com}}
  \country{} 
}

\begin{abstract}
As a telecom provider, our company has a critical mission to maintain telecom services even during power outages.
To accomplish the mission, it is essential to maintain the power of the telecom base stations.
Here we consider a solution where electric vehicles (EVs) directly supply power to base stations by traveling to their locations. The goal is to find EV routes that minimize both the total travel distance of all EVs and the number of downed base stations. In this paper, we formulate this routing problem as a new variant of the Electric Vehicle Routing Problem (EVRP) and propose a solver that combines a rule-based vehicle selector and a reinforcement learning (RL)-based node selector.
The rule of the vehicle selector ensures the exact environmental states when the selected EV starts to move.
In addition, the node selection by the RL model enables fast route generation, which is critical in emergencies.
We evaluate our solver on both synthetic datasets and real datasets.
The results show that our solver outperforms baselines in terms of the objective value and computation time. 
Moreover, we analyze the generalization and scalability of our solver, demonstrating the capability toward unseen settings and large-scale problems.
Check also our project page: \url{https://ntt-dkiku.github.io/rl-evrpeps}.
\end{abstract}


\newcommand{\BibTeX}{\rm B\kern-.05em{\sc i\kern-.025em b}\kern-.08em\TeX}

\begin{document}


\pagestyle{fancy}
\fancyhead{}


\maketitle 

\begin{figure}[tb] \centering
    \includegraphics[width=\linewidth]{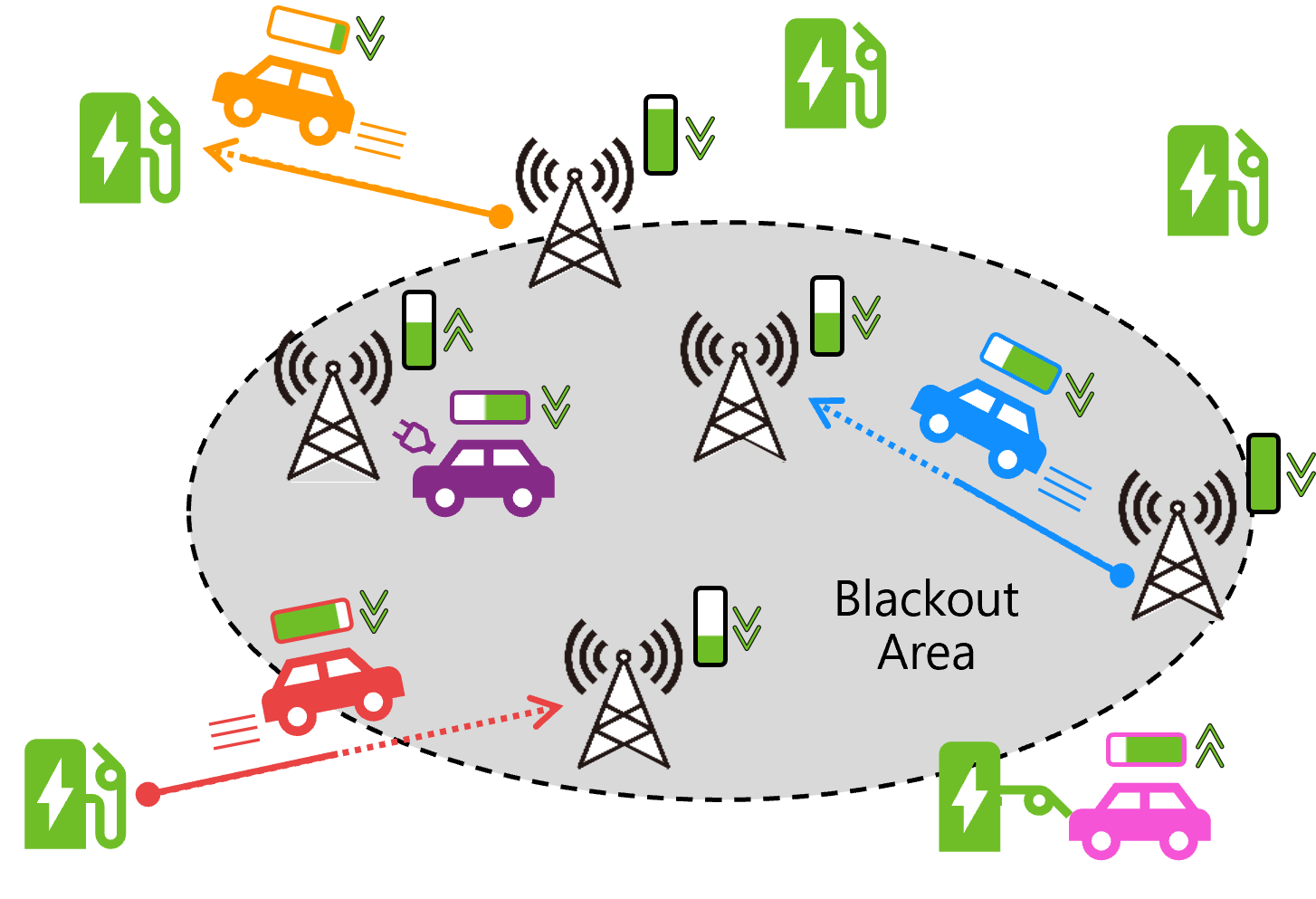}
    \caption{An illustration of telecom base station relief. EVs directly supply power to base stations by traveling to their locations. Each EV returns to a charge station located outside the blackout area before the EV power runs out.}
    \label{fig:teaser}
\end{figure}
\section{Introduction}
As the frequency of natural disasters increases \cite{disaster_report}, maintaining infrastructures during a disaster is becoming more critical. In particular, telecom services are one of the most important infrastructures to be always maintained as the Internet is nowadays a lifeline for people. As a telecom provider, our group company is committed to initiatives to keep providing telecom services even during a disaster.  
In this paper, we introduce an initiative to maintain the power of telecom base stations during power outages, which is one of the fundamental tasks of maintaining telecom services during a disaster with power outages. 
Each base station possesses a spare battery for short-time backup (e.g., around three hours), but an additional power supply from some external sources is required to maintain the base station batteries over a longer time.

Electric vehicles (EVs) are promising candidates for those external sources.
EVs can supply their power to objects using vehicle-to-everything (V2X),\footnote{V2X here refers to "electricity transfer" from EVs to everything, not "communication".} which has recently gained attention as a mobile power source for auxiliary services. 
In the context of emergency power supply, EVs are superior to conventional power-supply vehicles in the following: 
Companies can use EVs as company vehicles in normal times, thereby reducing waste of resources;
More EVs can be dispersed over a wider area because the cost per EV is less expensive.
Although the amount of electricity transported per EV in a single route is less than that of a power-supply vehicle, EVs are more suitable in our problem settings, given the need to power many dispersed base stations individually.

Recently, various approaches have emerged that apply vehicle-to-grid (V2G) or vehicle-to-home (V2H) technologies for emergency power supply \cite{EPS3, EPS2, EPS1}.
However, due to its high installation cost, introducing the grid system to numerous existing base stations is impracticable. 
Furthermore, one-to-one dispatch in \cite{EPS1} is inefficient in the base station relief where the number of EVs is smaller than the number of base stations to be rescued.
Hence, apart from these existing approaches, a routing approach where EVs maintain the power of the base stations by going around to power the base stations directly should also be considered.

In this paper, we formulate this routing problem as a variant of the Electric Vehicle Routing Problem (EVRP) and propose a solver that combines a rule-based vehicle selector and a reinforcement learning (RL)-based node selector.
Figure \ref{fig:teaser} illustrates our problem (we name it EVRP for Emergency Power Supply (EVRP-EPS)).
The goal is to find EV routes that minimize both the total travel distance of all EVs and the number of downed base stations in a situation where the base station batteries continuously deplete over time.
Compared to existing EVRPs, which consider battery consumption and recharging of EVs, EVRP-EPS additionally considers the battery discharge of EVs and some details mandatory in practice.
Our solver generates EV routes by alternately repeating the two components: 
The vehicle selector selects the EV that can make the next move the earliest;  
The node selector then determines the next destination of the selected EV with the centralized stochastic policy.
The rule of the vehicle selector ensures the exact environmental states when the selected EV starts to move. 
Furthermore, the RL-based node selection enables fast route generation, even in large-scale problems, which is critical during emergencies.

We evaluate our solver on synthetic and real datasets that involve actual locations and specifications of base stations, charge stations, and EVs. 
The results show that our solver outperforms baselines (two naive approaches and a conventional heuristic) in terms of both the object value and computation time. Moreover, we analyze the generalization and scalability of our solver, demonstrating the capability towards unseen settings and large-scale problems.

Contributions of this paper are organized as follows.
\begin{itemize}
    \item We are the first to formulate base station relief as a new EVRP that additionally considers the discharging of EVs.
    \item We propose a solver that combines a rule-based vehicle selector, an RL-based node selector, and a new objective function that considers both travel distance and the number of downed base stations.
    \item We demonstrate the effectiveness of our solver on real datasets that involve actual locations and specifications of base/charge stations and EVs.
\end{itemize}

\begin{table}[]
    \scriptsize
    \centering
    \caption{Notation.}
    \label{tab:notation}
    \begin{tabular}{l@{\extracolsep{5pt}}l@{\extracolsep{5pt}}l}
        \hline
        Notation & Description & Unit\\
        \hline
        $N_\textrm{bs}, N_\textrm{cs}, N, N_\textrm{ev}$ & The number of base stations, charge stations, nodes, and EVs & - \\
        & $N=N_\textrm{bs} + N_\textrm{cs}$& \\
        $i, j, k, n$ & The indices of base stations, charge stations, EVs, and nodes & -\\
        & $n = 1, ..., i, ..., N_\textrm{bs}, N_\textrm{bs}+1, ..., N_\textrm{bs}+j, ..., N$ & \\
        $\textrm{bs}_i$, $\textrm{cs}_j$, $\textrm{ev}_k$ & $i$-th base station, $j$-th charge station, and $k$-th EV & -\\
        $a_k$ & The indices of actions taken by the $k$-th EV & - \\
        $\Pi_k(a_k)$, $\pi_k(a_k)$ & The random variable and observed value of indices of nodes & - \\
                               & visited by $k$-th EV at $a$-th action& \\
        $t, T$ & Continuous time and time horizon& h\\
        \hline
        $B(\textrm{bs}_i,t)$ & $i$-th base station's battery at the time $t$ & kWh\\
        $C(\textrm{bs}_i)$ & $i$-th base station's power consumption & kWh/h\\
        $Q(\textrm{bs}_i)$ & $i$-th base station's capacity & kWh \\
        \hline
        $D(\textrm{cs}_j)$ & $j$-th charge station's discharge rate & kWh/h\\
        \hline
        $B(\textrm{ev}_k, t)$ & $k$-th EV's battery at the time $t$ & kWh\\
        $C(\textrm{ev}_k)$ & $k$-th EV's driving power consumption & kWh/km\\
        $Q(\textrm{ev}_k)$ & $k$-th EV's capacity & kWh \\
        $D(\textrm{ev}_k)$ & $k$-th EV's discharge rate & kWh/h\\
        \hline
    \end{tabular}
\end{table}

\section{Problem Setting}
In this section, we describe the EVRP-EPS, including the objective function, the action space, the battery fluctuation model of base stations and EVs, and other constraints.
Our notation is organized in Table \ref{tab:notation}.
\subsection{Objective Function}
Given a time horizon $T$ (e.g., expected blackout duration) and a set of base stations, charge stations, and EVs, the goal is to maintain as many base station batteries as possible during the time horizon while minimizing the total travel distance of all EVs. Formally, the objective function below is minimized.
\begin{equation}
    \label{eq:objective}
    \sum_k\sum_{a_k}^{A_k-1}\frac{\textrm{d}(\bm{x}_{\pi_k(a_k)}, \bm{x}_{\pi_k(a_k+1)})}{N_\textrm{ev}} + \alpha\frac{1}{T}\int_{t=0}^T\frac{\sum_{i}\textrm{I}(B(\textrm{bs}_{i}, {t}) = 0)}{N_\textrm{bs}}dt, 
\end{equation}
where $A_k$ the number of $k$-th EV’s actions, $\textrm{d}(\cdot, \cdot)$ is the distance between two points, $\bm{x}_n$ is the 2d coordinates of the $n$-th node, $\pi_k(a_k)$ is the index of node visited by $k$-th EV at $a$-th action, $\alpha$ is the positive weighting factor, and $\textrm{I}(\cdot)$ is the Boolean indicator function.
The first and second terms of Eq. (\ref{eq:objective}) correspond to the travel distance of all EVs and the number of downed base stations, respectively. 
\begin{figure}[tb] \centering
    \includegraphics[width=0.9\linewidth]{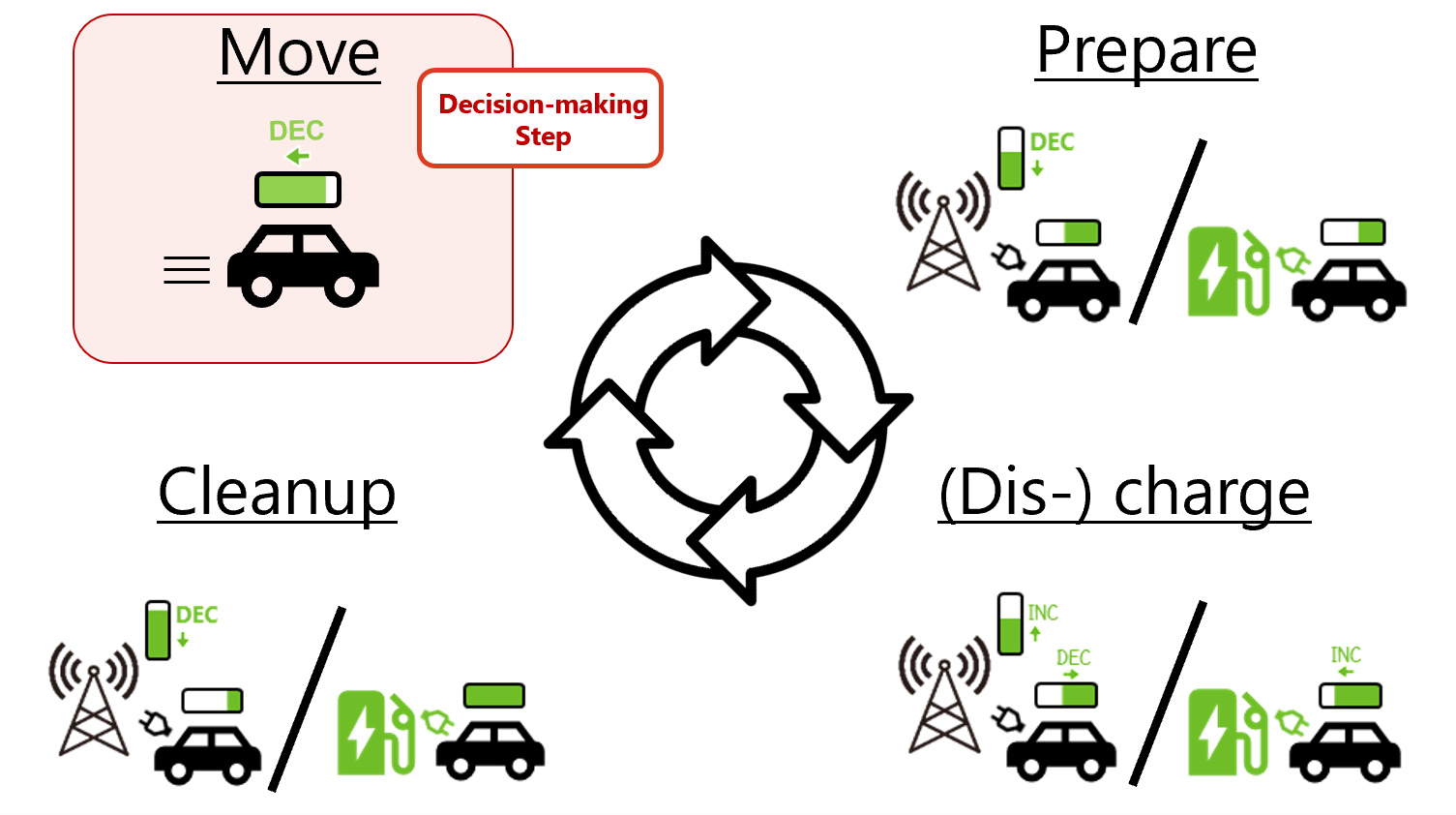}
    \caption{An illustration of the EV action cycle: \textit{move} $\to$ \textit{prepare} $\to$ \textit{(dis-) charge} $\to$ \textit{clean-up} $\to$ \textit{move} $\to$ $\cdots$}
    \label{fig:cycle}
\end{figure}

\begin{figure*}[tb] \centering
    \includegraphics[width=\textwidth,height=0.3\textwidth]{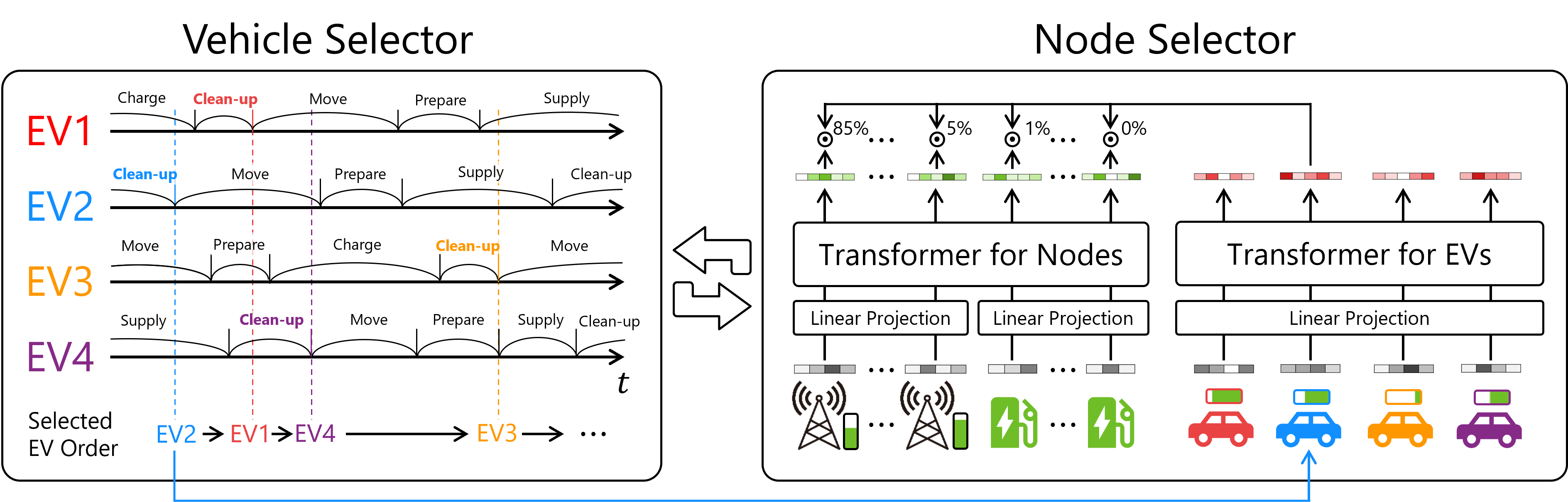}
    \caption{An overview of our solver. It first selects the EV that will be the earliest available for the next move. Then, it determines the next destination of the selected EV with a centralized stochastic policy, which is parameterized by a two-tower Transformer.} 
    \label{fig:method_overview}
\end{figure*}
\subsection{(Sub-) Actions and Action Space}
EVs/drivers cycle through an action: \textit{move} and three subsequent sub-actions: \textit{prepare}, \textit{(dis-) charge}, and \textit{clean-up} (Figure \ref{fig:cycle}).
First, each EV \textit{moves} to a node (i.e., base station/charge station).
After arriving at the node, it \textit{prepares} the equipment to (dis-) charge power.
Then it \textit{(dis-) charge} power at a base station/charge station.
Lastly, it \textit{clean-ups} the equipment and \textit{moves} to another node, and so on.
EVs keep this cycle during a given time horizon.


Here, the action space is solely about determining which node each EV moves to from the current node in the \textit{move} phase. 
The sub-actions (i.e., how long (much) EVs (dis-) charge power and how long preparation/cleanup takes) are incidentally determined by the environment states when the EV visits that node and constraints.


\subsection{Battery Fluctuation Model}
The battery of base stations/EVs linearly and continuously changes with time.
Let $\tau$ be continuous elapsed time, the battery of $i$-th base station changes depending on conditions, as follows:
\begin{equation}
    B(\textrm{bs}_{i}, {t+\tau}) = \left\{
        \begin{array}{ll}
            B(\textrm{bs}_{i}, t) + (D(\textrm{ev}_{k}) - C(\textrm{bs}_{i})) \tau & \textrm{if ev$_k$ discharges} \\[-0.5ex]
                                                                                   & \textrm{power at bs$_i$,} \\[1.1ex]
            B(\textrm{bs}_{i}, t) - C(\textrm{bs}_{i})\tau & \textrm{otherwise},
        \end{array}
    \right. 
\end{equation}
where $D(\textrm{ev}_{k}) - C(\textrm{bs}_{i}) > 0$ always holds, meaning that the amount of power supply from an EV to a base station always exceeds the power consumption of that base station.
During the $k$-th EV discharges power at the $i$-th base station, the battery of that base station increases in proportion to the elapsed time with the $k$-th EV discharge rate minus the $i$-th base station power consumption as a coefficient.
Otherwise (i.e., during no EV arrives at or an EV prepares/clean-ups at $i$-th base station), the battery of that base station decreases in proportion to the elapsed time with the power consumption as a constant.

On the other hand, the battery of $k$-th EV changes as follows:
\begin{equation}
    B(\textrm{ev}_{k}, t+\tau) = \left\{
        \begin{array}{ll}
            B(\textrm{ev}_{k}, t) - C(\textrm{ev}_{k})V\tau & \parbox{.3\linewidth}{if ev$_k$ moves} \\
            B(\textrm{ev}_{k}, t) - D(\textrm{ev}_{k}) \tau & \textrm{if ev$_k$ discharges power}\\
            B(\textrm{ev}_{k}, t) + D(\textrm{cs}_{j}) \tau & \textrm{if ev$_k$ charges at cs$_j$} \\
            B(\textrm{ev}_{k}, t) \               & \textrm{otherwise},
        \end{array}
    \right.
\end{equation}
where $V$ is the constant EV driving speed (km/h).
During $k$-th EV moves, its battery decreases proportionately to the travel distance $V\tau$ with the driving power consumption as a constant.
During $k$-th EV charges power from $j$-th charge station, its battery increases proportionately to the elapsed time with the discharge rate of that charge station as a constant. Otherwise (i.e., during $k$-th EV prepares/clean-ups), the EV battery does not change.
\subsection{Duration of EV (Sub-) Actions}
The duration of each EV (sub-) action $\tau'$ is deterministically computed as follows:  
\begin{equation}
    \label{eq:action_duration}
    \tau' = \left\{
        \begin{array}{ll}
            \frac{\textrm{d}\left(\bm{x}_{\pi_k(a)}, \bm{x}_{\pi_k(a+1)}\right)}{V}  & \textrm{if ev$_k$ moves,} \\[1.1ex]
            \frac{\beta Q(\textrm{ev}_k) - B(\textrm{ev}_k, t_\textrm{st})}{D(\textrm{cs}_j)} & \textrm{if ev$_k$ charges at cs$_j$,}\\[1.1ex]
            \min\left(\frac{\gamma Q(\textrm{bs}_i) - B(\textrm{bs}_i, t_\textrm{st})}{D(\textrm{ev}_k) - C(\textrm{bs}_i)}, \right. & \textrm{if ev$_k$ discharges power}\\[-1ex]
            & \textrm{at bs$_i$,}\\[-1ex]
            \ \ \ \ \ \ \ \ \left.\frac{B(\textrm{ev}_k, t_\textrm{st}) - Lim(\textrm{ev}_k, \textrm{bs}_i)}{D(\textrm{ev}_k)}\right) & \\[1.1ex]
            \Delta t_\textrm{pre}, \Delta t_\textrm{cln}  & \textrm{if ev$_k$ prepares/clean-ups,}\\
        \end{array}
    \right.
\end{equation}
where $t_\textrm{st}$ is the time at which each action starts, and $\Delta t_\textrm{pre/cln}$ is a pre-defined constant time for EVs preparing/cleanuping.
The duration of \textit{move} (travel time) is calculated assuming that the vehicle speed is constant. 
The duration of \textit{charge} is the time it takes for the EV to be charged up to $\beta\in[0, 1]\times100\%$ of its capacity.
The duration of \textit{discharge} is the time it takes for the base station battery to be charged up to $\gamma\in[0, 1]\times100\%$ of its capacity if the EV battery is sufficient. Otherwise, the time it takes for the EV to discharge until its discharge limit is reached. 
As the charge speed decreases after exceeding 80\% of the battery, we usually set $\beta$ and $\gamma$ to 0.8.
The $k$-th EV discharge limit at $i$-th base station $Lim(\textrm{ev}_k, \textrm{bs}_i)$ is calculated as,
\begin{equation}
    \label{eq:constraint1}
    Lim(\textrm{ev}_k, \textrm{bs}_i) = \max\left(\mu Q(\textrm{ev}_k), V\cdot\min_j\left(\textrm{d}(\bm{x}_{\textrm{bs}_i}, \bm{x}_{\textrm{cs}_j})\right)\right),
\end{equation}
where $\mu\in[0,1]$ is the discharge lower limit ratio to the EV capacity. 
Eq (\ref{eq:constraint1}) ensures that the EV stops discharging power when the discharge limit is reached, so that it is able to return to at least one depot without running out of its battery on the way. 
\subsection{Other Constraints}
The states of base stations and EVs always satisfy the following constraints:
\begin{align}
    &0 \leq B(\textrm{bs}_i, t) \leq Q(\textrm{bs}_i), \\
    \label{eq:constraints2}
    &C(\textrm{ev}_k)\cdot\left(\textrm{d}(\bm{x}^t, \bm{x}_n) + \min_j\left(\textrm{d}(\bm{x}_n, \bm{x}_{\textrm{cs}_j})\right)\right) \leq B(\textrm{ev}_k,t) \leq Q(\textrm{ev}_k),
\end{align}

The other unformulated constraints are summarized below.
\begin{itemize}
    \item In the \textit{move} phase, EVs do not choose base stations at which another EV stays, i.e., no more than 2 EVs stay at the same base station at the same time.
    \item If more than two EVs stay at the same charge station at the same time, EVs charge in the order in which they have arrived, with the $(q+1)$-th EV in the queue staying (doing nothing) until the $q$-th EV has finished charging.
\end{itemize}
\section{Methodology}
Here, we introduce our solver for the EVRP-EPS, which has been defined in the previous section.
Our solver generates EV routes by repeating the two components, vehicle selector and node selector (Figure \ref{fig:method_overview}).
In the vehicle selector, an EV is selected that finishes clean-up the soonest from the current time.
In the node selector, the next destination node of the selected EV is determined by a two-tower Transformer-based model, which parameterizes a stochastic policy.
After the selected EV starts to move to the determined destination node, another EV will be selected that finishes clean-up the soonest, and so on.
In the following, we describe the details of each component, how to train our solver, and how to generate the final routes.
\subsection{Vehicle Selector}
In the EV action cycle, EVs can move to the next destination node only after finishing the \textit{clean-up}.
Therefore, we employ a rule-based vehicle selector that always selects an EV that finishes \textit{clean-up} the soonest from the current time.
Formally, the index of the selected vehicle $\bar{k}$ is obtained by the Argumin function of the unmovable duration of all EVs:
\begin{equation}
    \bar{k} = \argmin_{k}{\left[U(\textrm{ev}_k, t)\right]},
\end{equation}
where $U(\textrm{ev}_k, t)$ is the unmovable duration of $k$-th EV at the current time $t$.
It is updated when the $k$-th EV starts moving as follows:
\begin{equation}
    U(\textrm{ev}_k, t) =  
    \left\{
        \begin{array}{ll}
            \tau'_\textrm{move} + \Delta t_\textrm{pre} + \tau'_\textrm{discharge} + \Delta t_\textrm{cln} & \textrm{If ev$_k$ moves to} \\[-1.0ex]
            & \textrm{a base station,} \\[1.1ex]
            \tau'_\textrm{move} + \Delta t_\textrm{pre} + \tau'_\textrm{charge} + \Delta t_\textrm{cln} & \textrm{If ev$_k$ moves to} \\[-1.0ex]
            & \textrm{a charge station,} \\
        \end{array}
    \right.
\end{equation}
where $\tau'_\textrm{move}$, $\tau'_\textrm{discharge}$, and $\tau'_\textrm{charge}$ are the duration of \textit{move}, \textit{discharge}, and \textit{charge} computed in Eq. (\ref{eq:action_duration}).
At other times it decreases with time as:
\begin{equation}
    U(\textrm{ev}_k, t+\tau) = U(\textrm{ev}_k, t) - \tau.
\end{equation}
\subsection{Node Selector}
Here, we aim to obtain a centralized stochastic policy that outputs the probability of visiting each node, given the global environment state and the index of a selected EV. In this paper, we parameterize this stochastic policy with a two-tower Transformer-based model (the r.h.s of Figure \ref{fig:method_overview}).

Our model first computes the initial embeddings of nodes and EVs with linear projections:
\begin{align}
    \bm{h}_{\textrm{node}_n}^{(0)} &= \left\{
        \begin{array}{ll}
             W_\textrm{bs}\bm{s}({\textrm{bs}_n, t_{a_{\bar{k}}}}) + \bm{b}_\textrm{bs}&  \textrm{if}\ n \in \{1,\dots,N_\textrm{bs}\}, \\
             W_\textrm{cs}\bm{s}({\textrm{cs}_{n-N_{\textrm{bs}}}, t_{a_{\bar{k}}}}) + \bm{b}_\textrm{cs}& \textrm{if}\ n \in \{N_\textrm{bs}+1,\dots, N\}, 
        \end{array}
    \right.\\
    \bm{h}_{\textrm{ev}_k}^{(0)} &= W_\textrm{ev}\bm{s}(\textrm{ev}_k, t_{a_{\bar{k}}}) + \bm{b}_\textrm{ev},
\end{align}
where $W_\textrm{bs}, W_\textrm{cs}, W_\textrm{ev}, \bm{b}_\textrm{bs}, \bm{b}_\textrm{cs}, \bm{b}_\textrm{ev}$ are trainable projection matrices and biases for base stations, charge stations, and EVs, respectively. 
$t_{a_{\bar{k}}}$ is the time when the selected EV takes the $a$-th action (right after finishing \textit{clean-up}).
The input base-station state $\bm{s}(\textrm{bs}_n, t)$ includes 2d coordinates, capacity, consumption rate, remaining battery at $t$, and expected time to be downed (i.e., the feature dimension $d=6$). 
The charge-station state $\bm{s}(\textrm{cs}_{n-N_\textrm{bs}}, t)$ includes 2d coordinates, discharge rate, and boolean of whether being visited by an EV (i.e., $d=3$). The EV state $\bm{s}(\textrm{ev}_k, t)$ includes 2d coordinates, a boolean of whether being at a charge station, current phase in the EV cycle, duration of each phase, unmovable duration, capacity, and remaining battery at $t$ (i.e., $d=12$).

It then produces the final embeddings of nodes and EVs by stacking $L$ of the following Transformer encoders \cite{Transformer}.
\begin{align}
    \bm{h}^{(l)}_{\textrm{node}_n} &= \textsc{Xfmr}_{\textrm{node}_n}^{(l)}\left(\bm{h}^{(l-1)}_{\textrm{node}_1},\dots,\bm{h}^{(l-1)}_{\textrm{node}_N}\right), \\
    \bm{h}^{(l)}_{\textrm{ev}_k} &= \textsc{Xfmr}_{\textrm{ev}_k}^{(l)}\left(\bm{h}^{(l-1)}_{\textrm{ev}_1},\dots,\bm{h}^{(l-1)}_{\textrm{ev}_K}\right),
\end{align}
where $\textsc{Xfmr}_{\textrm{node}_n}^{(l)}, \textsc{Xfmr}_{\textrm{ev}_k}^{(l)}$ are the $l$-th Transformer encoders for nodes and EVs, of which subscript indicate the output element. Note that no positional encoding is used here as nodes and EVs are permutation-invariant.
\par
Finally, the (conditional) probability of visiting each node is computed from the scaled dot-product attention between the final embeddings of nodes and the selected EVs, as follows.
\begin{equation}
    u(k, n) = \left\{
    \begin{array}{ll}
        \eta\cdot\tanh\left(\frac{\bm{q}^{\top}_{\textrm{ev}_k}\bm{k}_{\textrm{node}_n}}{\sqrt{d_\textrm{k}}}\right)  & \textrm{if node$_n$ is visitable,}  \\
        -\infty & \textrm{otherwise,} 
    \end{array} \right.
\end{equation}
\begin{equation}
    p_\theta\left(\Pi_{{\bar{k}}}(a_{\bar{k}})=n|{\mathcal{S}_{t_{a_{\bar{k}}}}},\bar{k}\right) = \frac{e^{u(\bar{k}, n)}}{\sum_m e^{u(\bar{k}, m)}},
\end{equation}
where $\eta (=10)$ is the clipping width, the query $\bm{q}_{\textrm{ev}_k} = W^Q\bm{h}_{\textrm{ev}_k}^{(L)}$, the key $\bm{k}_{\textrm{node}_n} = W^K\bm{h}_{\textrm{node}_n}^{(L)}$, $d_\textrm{k}$ is the dimension of the key, $W^Q, W^K$ are trainable projection matrices, and $\mathcal{S}_t$ is the global state that consists of the states of all base stations, charge stations, and EVs at $t$. 
Regarding the conditional branch, a base station is visitable if the selected EV can return to a charging station without running out of battery after visiting that base station, and no other EVs visit it.
A charge station is visitable if the selected EVs can reach it from the current node without running out of battery on the way.
\subsection{Training}
We train the stochastic policy with REINFORCE \cite{REINFORCE}, a policy gradient algorithm. As the baseline, we employ the greedy rollout baseline similar to \cite{AM}. The gradient $g$ is computed as follows, 
\begin{equation}
    g = \sum_b \left(\mathcal{L}(\bm{\pi}_b) - \mathcal{L}(\bm{\pi}_b^{GR})\right) \nabla_{\theta}\log p_{\theta}(\bm{\pi}_b),
\end{equation}
where $\mathcal{L}$ is the cost function of Eq. \eqref{eq:objective}, $\bm{\pi}_b$ is the route generated by sampling rollout on the current training policy, $\bm{\pi}^{\textrm{GR}}_b$ is the route generated by greedy rollout on the baseline policy.
The baseline policy is replaced with the current training policy at the end of each epoch only if the improvement between the two policies is more than 5 \% in a paired t-test on validation instances. 
The subscript $b$ indicates the route is for $b$-th instance of a batch. $p_{\theta}(\bm{\pi})$ is the probability of generating the route $\bm{\pi}$, which is factorized as:
\begin{equation}
    p_{\theta}(\bm{\pi}) = \prod_k\prod_{a_k}{p_{\theta}(\pi_k(a_k)|\mathcal{S}_{t_{a_k}}, k)}.
\end{equation}

\subsection{Decoding (Route Generation)}
We can generate the final route by sampling with the trained stochastic policy.
The sampling technique includes greedy decoding, sampling decoding, beam search, and Monte Carlo tree search.
In this paper, we employ the sampling decoding, which is simple yet effective.
It first samples around 1k - 10k routes with sampling rollout, then selects the route that minimizes the cost function the most as the final route.


\section{Experimental set-ups}
\subsection{Datasets}
We evaluate our solver on synthetic datasets and real datasets. 
The synthetic datasets are divided into training, validation, and evaluation splits. 
We train our solver on the training split and select the weights of the epoch where the cost function is the lowest on the validation split.
Note that we use greedy decoding in this validation phase.
We then evaluate it on the evaluation split and real datasets.

\subsubsection*{\textbf{The real datasets (\textsc{Real-ev-6, 12})}} 
They include actual data of two different regions\footnote{\textbf{The names of the regions and specific values are omitted for privacy reasons.}}: the areas of \textsc{Real-ev-6} and \textsc{Real-ev-12} are around 46 and 66 km$^2$, respectively. 
Table \ref{tab:statics} organizes the statistics of these datasets.
The locations and specifications of nodes and EVs are actual values.
Base stations are distributed in the center and charge stations surround it, emulating the actual situation.
We assume that the initial batteries of base stations and EVs are 80 \% of their capacities.

\subsubsection*{\textbf{The synthetic datasets (\textsc{Syn-ev-6, 12})}} 
They are generated based on the distribution of the real datasets. Table \ref{tab:statics} organizes their statistics.
The locations of nodes are uniformly sampled within the square area $[0, 100]^2$ km$^2$. 
The specific values of base stations, charge stations, and EVs are sampled from the distribution of real datasets.
The initial batteries of base stations are randomly sampled within 50 - 100 \% of the capacity, and those of EVs are set to 80 \% of the capacity. 
\textsc{Syn-ev-6S} is used only for scalability and generalization tests, which is discussed later.

\subsubsection*{\textbf{Parameters of base stations, charge stations, and EVs}}
We here consider the Nissan Leaf e+ series, of which the capacity $Q(\textrm{ev}_k)$ = 60 kWh, driving power consumption $C(\textrm{ev}_k)$ = 161 Wh/km, discharge lower limit $\mu$ = 0.1. The discharge rate of EVs is set as $D(\textrm{ev}_k)$ = 10 kWh/h. The constant vehicle speed is set as $V$ = 41 km/h. 
The charge upper limit is set to 80\%, i.e., $\beta,\gamma$=0.8.
The discharge rate of charge stations is either 3 (normal charging) or 50 kWh/h (rapid charging), i.e., $D(\textrm{cs}_j)$ = 3, 50 kWh/h.
Prepare/cleanup duration (i.e., $\Delta t_\textrm{pre/cln}$) at base stations and charge stations are set to 0.5h and 10m, respectively.
Each base station has its own capacity and consumption rate, details of which are omitted.
The distance between nodes is approximated by Euclidean distance.
\begin{table}[tb]
    \caption{The statistics and use of the synthetic datasets and real datasets.}
    \begin{center}
    \begin{tabular}{l @{\extracolsep{4pt}}cccc@{\extracolsep{4pt}}l@{}}
        \hline
        Name & $N_\textrm{ev}$ & $N_\textrm{bs}$ & $N_\textrm{cs}$ & Area (km$^2$) & Use\\
        \hline
        \textsc{Real-ev-6}  &6   &33   &62  & 46 & Eval\\
        \textsc{Real-ev-12} &12  &46   &117 &66 & Eval\\
        \hline
        \textsc{Syn-ev-6}   &6   &50  &12   &100 & Train/Valid/Eval\\
        \textsc{Syn-ev-12}  &12  &50  &12   &100 & Train/Valid/Eval\\
        \hline
        \textsc{Syn-ev-12S} &12  &25  &12   &100 & Train/Valid\\
        \hline
    \end{tabular}
    \label{tab:statics}
    \end{center}
\end{table}

\subsection{Baselines}
We compare our solver with the following solvers.

\subsubsection*{\textbf{Random Node Selector (\textsc{Rand})}}
It replaces the node selector of our solver with a random selection. Given a selected EV, it randomly selects a visitable node as the next destination of that EV. The final route is the best route among 12.8k sampled routes, similar to the sampling decoding of our solver.

\subsubsection*{\textbf{Greedy Node Selector (\textsc{Greed})}}
It replaces the node selector of our solver with a greedy selection. Given a selected EV, it selects a visitable base station with the lowest current battery as the next destination of that EV. If there is no visitable base station, it selects the nearest charge station from the current position.

\subsubsection*{\textbf{Constraint Programming on Time-Space Network (\textsc{Cp-Tsn})}}
The Time-Space Network (TSN) \cite{TSN} duplicates nodes along the discrete-time direction and connects two nodes at different discrete times with an arc.
The Constraint Programming (CP) \cite{CP} then finds the optimal chain of the arcs, which corresponds to an EV route.
As the computation of the vanilla model is too expensive to derive a solution within the time limit, we first assign nodes grouped by a balanced k-means to each EV and apply \textsc{Cp-Tsn} to each of the divided problems separately. Note that solving divided problems restricts the search space and may prevent deriving the original problem's optimal solution (settling for a near-optimal solution).
The OR-Tools CP-SAT solver\footnote{\url{https://developers.google.com/optimization/cp/cp_solver}} is used here. 
We report the results of two different time resolutions, 1h and 0.5h, denoted by $\textsc{Cp-Tsn}_{1.0}$ and $\textsc{Cp-Tsn}_{0.5}$, respectively. 


\subsection{Hyperparameters and Devices}
Trainable parameters are initialized with $\mathcal{U}(-1/\sqrt{d}, 1/\sqrt{d})$, $d$ is the input dimension. 
Training, validation, and evaluation splits for the synthetic datasets have 1.28M, 1000, and 100 instances, respectively. 
We evaluate a solver trained with the same number of EVs and $T=12$ as the evaluation data.
The maximum epoch is 100, the mini-batch size is 256, and the constant learning rate is $10^{-4}$.
The number of Transformer encoder layers is set as $L=2$, the dimension in hidden layers $H=128$, and the number of heads $M=8$.
The weighting factor is set as $\alpha=100$.
The random seed is set to 1234 unless otherwise stated.
A single GPU (RTX A6000: 48G) and two CPUs (Xeon Platinum 8380 (2.30GHz)) are used in the experiments. 
\section{Experimental results}
\begin{table*}[tb]
    \caption{Our solver (Ours) v.s. baselines. Metrics are averaged travel distance per EV (dist (km)), the time average of \# downed base stations (down), the objective value (obj), and total computation time (time). For all values, smaller is better.}
    \small 
    \begin{center}
    \begin{tabular}{l cccc c cccc c cccc c cccc}
        \hline
        &\multicolumn{9}{c}{\textsc{Syn-ev6} (100 samples)} & &\multicolumn{9}{c}{\textsc{Syn-ev12} (100 samples)} \\
        \cline{2-10} \cline{12-20}
        &\multicolumn{4}{c}{$T$ = 12h} &&\multicolumn{4}{c}{$T$ = 24h}  &&\multicolumn{4}{c}{$T$ = 12h} &&\multicolumn{4}{c}{$T$ = 24h}\\
        \cline{2-5} \cline{7-10} \cline{12-15} \cline{17-20}
        Model & dist & down & obj & time && dist & down & obj & time && dist & down & obj & time && dist & down & obj & time \\
        \hline
        w/o EVs &- &20.1 &- &- &&- &33.3 &- &- &&- &20.1 &- &- &&- &33.3 &- &-\\
        \textsc{Greed}  &189 &17.7 &37.2 &1s   &&312 & 29.8 &62.8 &1s &&192 &15.5 &32.9 &1s   &&317 &26.6 &56.4 &2s\\
        \textsc{Rand} (S=12800) &142 &15.3 &\textbf{32.0} &1m &&263 &26.5  &\textbf{55.5} &2m &&157 &11.8 &\textbf{25.1} &3m &&283 &21.1 &\textbf{45.1} &6m\\
        \hline
        Ours (G)       &81 &12.5 &25.9 &1s          &&125 &25.2 &51.7 &1s            &&78 &6.63 &14.0 &1s            &&131 &18.2 &37.8 &2s\\
        Ours (S=1280)  &78 &12.4 &25.5 &51s         &&125 &24.6 &50.4 &1m            &&75 &6.41 &13.6 &1m            &&129 &17.4 &36.0 &3m\\
        Ours (S=12800) &78 &12.4 &\textbf{25.5} &9m &&126 &24.5 &\textbf{50.2} &15m  &&74 &6.38 &\textbf{13.5} &15m  &&129 &17.2 &\textbf{35.7} &28m\\
        \hline
        \hline
        &\multicolumn{9}{c}{\textsc{Real-ev6} (1 sample)} & &\multicolumn{9}{c}{\textsc{Real-ev12} (1 sample)} \\
        \cline{2-10} \cline{12-20}
        &\multicolumn{4}{c}{$T$ = 12h} &&\multicolumn{4}{c}{$T$ = 24h}  &&\multicolumn{4}{c}{$T$ = 12h} &&\multicolumn{4}{c}{$T$ = 24h}\\
        \cline{2-5} \cline{7-10} \cline{12-15} \cline{17-20}
        Model & dist & down & obj & time && dist & down & obj & time && dist & down & obj & time && dist & down & obj & time \\
        \hline
        w/o EVs &- &8.25 &- &- &&- & 17.5 &- &- &&- &11.3 &- &- &&- &25.9 &- &-\\
        \textsc{Greed} &30 &5.87 &19.4 &1s &&37 &13.5 &42.8 &1s &&88 &9.34 &21.7 &1s &&107 &20.9 &47.1 &1s\\
        \textsc{Rand} (S=12800) &32 &3.74 &\textbf{13.1} &1s &&63 &10.0 &\textbf{33.8} &2s &&120 &4.98 &\textbf{12.8} &5s &&215 &14.2 &\textbf{34.2} &8s\\
        \hline
        \textsc{Cp-Tsn}$_{1.0}$ &16 &5.62 &17.9 &1m            &&24 &12.6 &39.37         &30m &&47 &4.54 &10.62 &22s        &&98 &14.2 &\textbf{32.44} &30m\\
        \textsc{Cp-Tsn}$_{0.5}$ &16 &3.40 &\textbf{11.2} &11m  &&31 &11.6 &\textbf{36.7} &30m &&53 &3.32 &\textbf{8.06} &1m &&105 &15.0 &34.25 &30m\\
        \hline
        Ours (G)       &17 &2.53 &8.61 &1s          &&26 &8.46 &27.1 &1s           &&59 &1.06 &3.26 &1s           &&99  &11.1 &25.6 &1s\\
        Ours (S=1280)  &21 &2.03 &7.27 &1s          &&35 &8.37 &27.3 &2s           &&61 &0.65 &2.41 &4s           &&101 &10.8 &25.2 &6s\\
        Ours (S=12800) &19 &1.96 &\textbf{6.97} &9s &&33 &8.15 &\textbf{26.5} &15s &&61 &0.48 &\textbf{2.03} &28s &&100 &10.6 &\textbf{24.7} &47s\\
        \hline
    \end{tabular}
    \label{tab:results}
    \end{center}
\end{table*}
\begin{figure*}[tb] \centering
    \includegraphics[width=0.9718\textwidth]{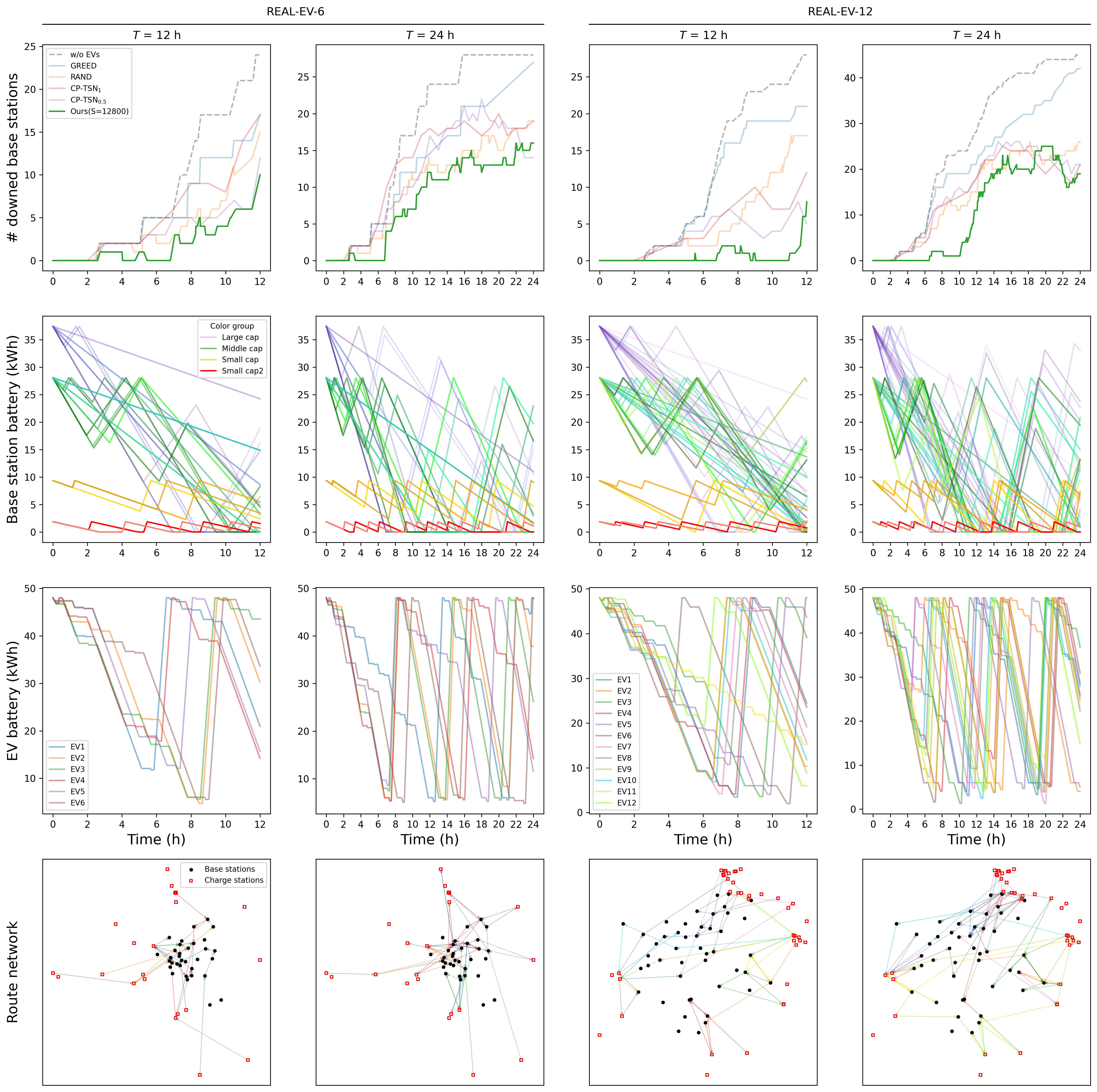}
    \caption{Four visual results for each real dataset: The time variation in the number of downed base stations (1st row),  the time variation in base station/EV batteries (2nd-3rd rows), and the route network (4th row). Each column corresponds to each of the real datasets and includes those corresponding results. }
    \label{fig:vis_routes}
\end{figure*}

\subsection{Performance Comparison with Baselines}
Table \ref{tab:results} shows the comparison results on the synthetic datasets and real datasets with $T$=12, 24h.
Metrics are averaged travel distance per EV (\texttt{dist} (km)), the time average of the number of downed base stations (\texttt{down}), the objective value (\texttt{obj}), and total computation time (\texttt{time}). 
Here, we set the time limit to 30 minutes, which is the requirement from the site workers.
When calculations exceed the time limit, we report the intermediate results obtained at that point.
In this comparison, \texttt{obj} and \texttt{time} are the primary metrics, with \texttt{dist} and \texttt{down} considered as reference information.\footnote{This is because \texttt{obj} already considers the trade-offs between \texttt{dist} and \texttt{down}.}

\subsubsection*{\textbf{Our solver (Ours)}}
Ours consistently outperforms baselines on all datasets in terms of both \texttt{obj} and \texttt{time}.
Notably, ours provides the best \texttt{obj} even on real datasets, where the distribution of nodes significantly differs from that of the training datasets (i.e., uniform distribution v.s. two clusters: one located at the center and the other surrounding the perimeter). These results demonstrate that ours has the generalization ability for node distribution and can be applied to any region.
In comparison between the decoding types of ours, the sampling decoding (Ours (S = \# samples)) outperforms the greedy decoding (Ours (G)). 
Furthermore, We observe an improvement in \texttt{obj} as the sample size increases.
The computation time also increases at that time, but there is enough margin for the time limit.

\subsubsection*{\textbf{\textsc{Greed} \& \textsc{Rand}}}
\textsc{Greed} always provides the worst \texttt{obj}, indicating that a simple rule-based routing is ineffective in EVRP-EPS.
By contrast, \textsc{Rand} provides relatively reasonable \texttt{obj}.
As these naive solvers do not consider minimizing travel distance, \texttt{dist} tends to be longer than other solvers.
The comparison between \textsc{Rand} and ours in the same sample size demonstrates that ours effectively samples better routes, guided by the trained stochastic policy.

\subsubsection*{\textbf{\textsc{Cp-Tsn}}}
It provides the second-best \texttt{obj} on real datasets.
We observe an improvement in \texttt{obj} as the time resolution increases.
This is because higher time resolution approximates continuous time more closely, thereby reducing the redundant waiting periods that arise between discrete time ticks and the completion of an action.
However, as the search space grows exponentially with the increase in time resolution,
we should keep it at a reasonable value (e.g., 0.5-1h) to find a reasonable solution within the time limit.
Furthermore, \textsc{Cp-Tsn} solves sub-problems that are not guaranteed to recover the original problem fully, and therefore that may be why its performance is inferior to ours. Overall, the poor scalability of this solver is a major limitation to derive a reasonable route within a constrained timeframe.

\subsubsection*{\textbf{Visualization}}
Figure \ref{fig:vis_routes} shows four types of visual results on the real datasets: The time variation in the number of downed base stations, the time variation in base station/EV batteries, and the route network formed by EVs' travel trajectory.

The time variation in the number of downed base stations provides more fine-grained perspectives of \texttt{down}, which corresponds to (the area bounded below the curve and by the $x$-axis) / (time horizon $T$) here.
Ours maintains the lowest number of downed base stations compared to other solvers for most of the time horizon.
The reduction by ours is greater in \textsc{Real-EV12}, where the ratio of base stations to EVs is smaller than in \textsc{Real-EV6} (3.8 v.s. 5.5). This indicates ours works more effectively when the theoretical maximum number of base stations that can be maintained is greater.
Regarding the increase in that number towards the end of the time horizon, it is preferable to a mid-period rise, as we here assume power outages are resolved immediately after the end of the time horizon.

The time variation in base station batteries shows that there are four distinct base station capacities, and small-capacity base stations are visited many times in a short span of time. 
We also observe that EVs do not visit base stations that will not be down during the time horizon, even without a power supply from EVs.

In the time variation in EV batteries, we observe two distinct discharge patterns: one repeats a small amount of discharge numerous times, while the other performs a relatively substantial discharge around two to five times.
The former can be associated with discharging to small-capacity base stations that require frequent charging, while the latter corresponds to discharging to other middle to large-capacity base stations.
Despite the absence of explicitly set roles, these two roles arise because once an EV visits a small-capacity base station, it becomes rational for it to continue visiting other small-capacity stations:
The shorter duration of discharging at small-capacity base stations allows the EV there to move on to the next node sooner than other EVs;
It is a natural choice for that EV to visit other small-capacity base stations at that time as they are likely down within a short time;
As a result of this repeated, that EV visits small-capacity base stations many times while other EVs continue to discharge at one base station. 

The visualization of route networks shows two node distinct distributions of \textsc{Real-EV6,12} and the travel trajectory of each EV with different colors.
Intuitively, there are no obviously inefficient visits, and all routes are considered reasonable.
The route networks with small-capacity base stations as hubs also confirm that these stations are visited frequently.

\subsection{Computational Scalability Test}
One of the biggest reasons why we employ a RL-based approach is its fast computation for large-scale problems.
To confirm this, we here evaluate the computational scalability of our solver: How much the computation time increases with respect to the increase in the time horizon and the number of base stations and EVs. 
We prepare additional evaluation datasets, which are organized in Table \ref{tab:scale_gen_test}. 
They are generated with the random seed 0 and the sample size 100.
The training dataset for  \textsc{Syn-Th-*} and \textsc{Syn-Bs-*} is \textsc{Syn-ev-12S}, while for \textsc{Syn-Ev-*}, it is \textsc{Syn-ev-6}.
The top row of Figure \ref{fig:scale_general_test} shows the computation time in each dataset. Assuming actual situations, the computation time here refers to the average time per sample when deriving one route at a time (i.e., the batch size is 1).
We observe different behaviors in the increase of each parameter.
The computation time increases approximately in proportion to the time horizons.
This is because the length of the action sequence $A_k$ increases in proportion to the time horizons. 
With respect to the increase in the number of base stations, the computation time is kept below the linear increase.
This indicates that batch parallelization with GPUs and other optimizations suppress, on this (practical) scale, the quadratic time complexity with respect to the number of nodes.
The computation time with respect to the number of EVs increases rapidly compared to other parameters.
It appears that the length of the action sequence increases proportionally with the number of EVs, 
but this is not actually the case. 
In fact, more EVs lead to frequent small discharges of each EV, resulting in more computation time (i.e., actions) than expected.
However, there is still enough margin for the time limit (30 minutes).
Overall, these results demonstrate that our solver generates a route at every conceivable practical scale with enough speed for the time limit.

\begin{table}[tb]
    \caption{The statistics of additional evaluation datasets (100 samples/dataset) for scalability and generalization tests. Note that \textsc{Syn-Th-*} are all the same dataset, where only the time horizon during inference is set differently.}
    \begin{center}
    \footnotesize
    \begin{tabular}{l cccc}
        \hline
        Name & $N_\textrm{ev}$ & $N_\textrm{bs}$ & $N_\textrm{cs}$ & $T$\\
        \hline
        \textsc{Syn-Th-12,24,36,48} &12  &25   &12 & \{12, 24, 36, 48\}\\
        \textsc{Syn-Bs-25,50,75,100} &12  &\{25, 50, 75, 100\} &12 &12\\
        \textsc{Syn-Ev-6,12,18,24} &\{6, 12, 18, 24\}  &50  &$N_\textrm{ev}$ &12\\
        \hline
    \end{tabular}
    \label{tab:scale_gen_test}
    \end{center}
\end{table}

%

\subsection{Generalization Test}
\begin{figure}[tb] \centering
    \includegraphics[width=\linewidth, height=0.6\linewidth]{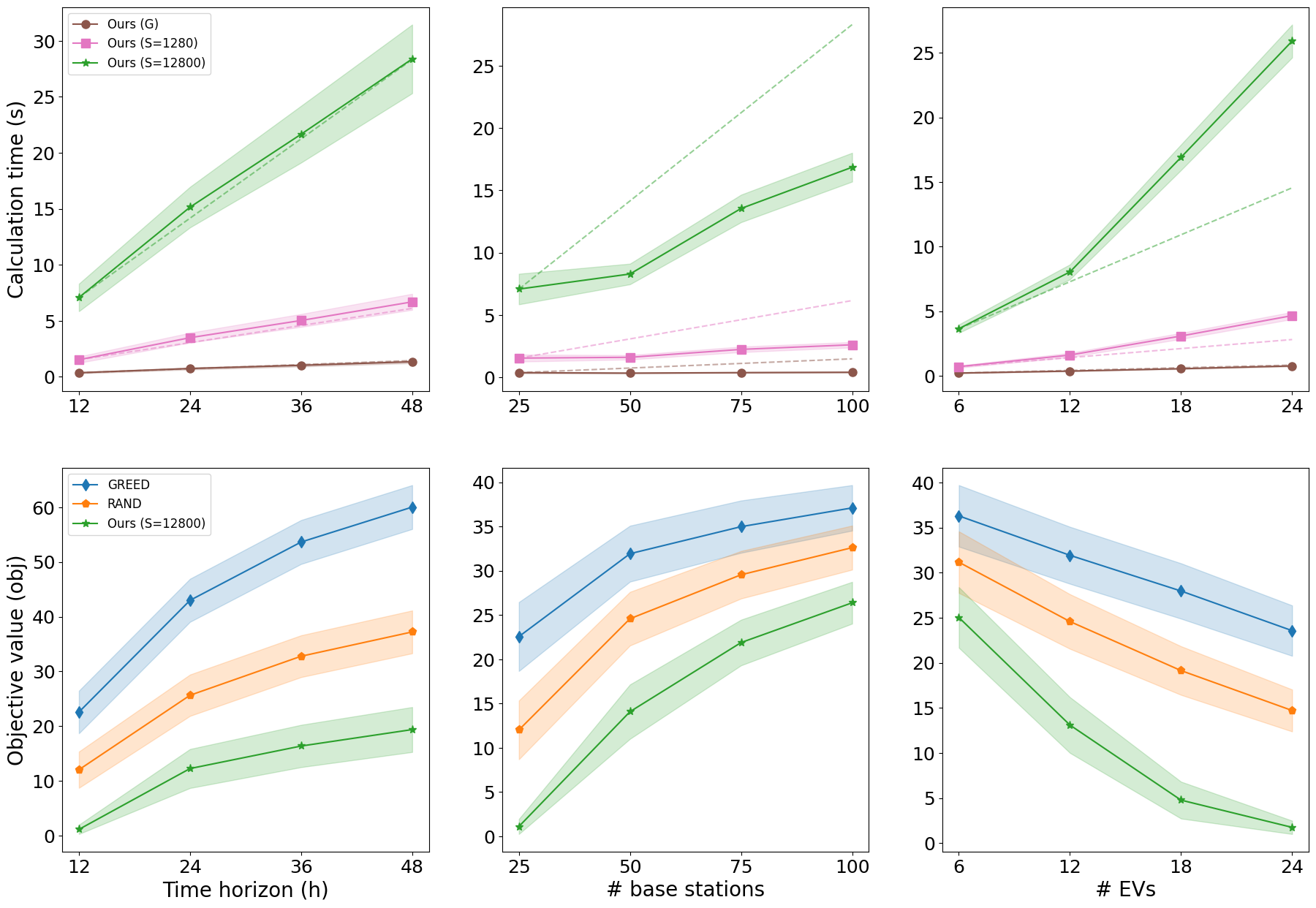}
    \caption{Results of scalability test (top row) and generalization test (bottom row). The points represent average values. The filled area indicates the standard deviation. The dotted line represents the linear increase from the starting point.} 
    \label{fig:scale_general_test}
\end{figure}

The time horizon and the number of EVs/nodes vary depending on the situation.
However, it is impracticable to train our solver for each case.
Therefore, our solver requires the generalization ability for unseen settings.
We here use the same datasets and our solvers as the scalability test.
The bottom row of Figure \ref{fig:scale_general_test} shows the objective values in each dataset, which are evaluated at the same time as the computation time in the scalability test.
We report the comparison between our solver and naive approaches (\textsc{Greed} and \textsc{Rand}).\footnote{We here exclude \textsc{Cp-Tsn}$_{*}$ because it exceeds the index limit of OR-Tools (i.e., outputs \texttt{IndexError: list index out of range}) in the middle of increasing the parameters.}
Our solver consistently outperforms the naive approaches at every parameter setting, where our solver is trained on the most-left parameter setting in each graph.
Notably, we observe that our solver effectively increases the improvement ratio in the unseen settings of 18 and 24 EVs.
These results demonstrate that our solver possesses reasonable generalization capabilities for unseen settings.
Besides, they confirm that in the previous scalability test, our solver provided not only fast computation but also adequate performance.


\section{Related Work}
\subsection{Vehicle Routing Problem}
The Vehicle Routing Problem (VRP) is a combinatorial optimization problem that aims to find the optimal routes for a fleet of vehicles to visit a given set of destinations.
Since Dantzig and Ramser \cite{vrp_initial} introduced its concept, numerous variants of VRPs have been proposed along with their solvers.
We will discuss related variant VRPs (i.e., EVRP) later, and here we focus on the types of solvers.

There exists a trade-off between the route quality and computation time, 
therefore, the appropriate solver must be chosen based on the situation.
Exact solvers such as the Branch-and-Cut (-and-Price) algorithm \cite{BC,exact_solver} provide the strictly optimal routes, while they have difficulty solving large-scale VRPs quickly.
On the other hand, heuristics \cite{LK, LKH, HGS, GA} provide near-optimal routes within a reasonable time, even in large-scale VRPs.
Recently, neural network-based solvers have been explored \cite{PtrNet, bello, nazari, ResGatedGCN, AM}, which provide faster computation than heuristics.
In addition, they can generate a heuristic optimized to a specific VRP from the VRP data, thereby overcoming the difficulty such that constructing specialized hueristics requires deep domain knowledge.
In this paper, we employ an NN-based solver due to the scale and time limit of EVRP-EPS.
The lack of specialized heuristics for EVRP-EPS also justifies our choice.

Among the NN-based solvers, Transformer-based multi-agent RL (TMARL) solvers are most relevant to our solver.
Several TMARL solvers have been proposed for multi-vehicle settings \cite{MARL_VRP1,MARL_VRP2,MARL_VRP3,MARL_VRP4}. 
However, they are not designed for our situation, where node property (i.e., battery) changes continuously over time, and these changes are influenced by the timing at which the node is visited (multiple visits per node are possible).
Our solver addresses this by introducing the rule-based vehicle selector, which allows it to follow the exact global state in our complicated situation. 
Decentralized learning is another choice to address this, but we employ centralized learning with the rule-based guide because we here assume the use of our centralized system that monitors the status of EVs and nodes.
Regarding the node selector, we propose the simplest two-tower model among existing TMARL solvers, of which performance remains sufficient.

\subsection{Electric Vehicle Routing Problem}
The Electric VRP (EVRP) is a VRP that additionally considers EV battery consumption by travel and recharging of EV battery.
Anticipating the introduction of EVs in logistics, various works related to EVRP have appeared lately (refer to \cite{EVRP_review1, EVRP_review2} for more comprehensive reviews on this topic).
Here, we focus on the problem formulation.
There are many variants of EVRP, where the conventional VRP vehicles are replaced with EVs, and their characteristics are considered.  
They include the Electric Traveling Salesman Problem with Time Windows (ETSPTW) \cite{ETSPTW}, EVRP \cite{EVRP}, EVRP with Time Windows (EVRPTW) \cite{EVRPTW,EVRPTW2,EVRPTW3}, and Electric Pickup and Delivery Problem with Time Windows (EPDPTW) \cite{EPDPTW}.
From more detailed perspective, battery consumption estimation according to the EV's load has also been studied to make the problem more realistic: linear deterministic \cite{linear_consump1, linear_consump2, linear_consump3}, non-linear deterministic \cite{nonlinear_consump1, nonlinear_consump2}, and stochastic estimation \cite{stat_consump1, stat_consump2, stat_consump3, stat_consump4}. 
Our formulation, EVRP-EPS, can be considered a variant of the Continuous Inventory Routing Problem (CIRP) \cite{IRP}.
However, it has not yet been studied in the context of EVs for emergency power supply, where the inventory corresponds to the base station battery.
Regarding battery consumption estimation in EVRP-EPS, linear models are used, but we found them sufficient based on actual driving data in our setting.

\subsection{EVs for Emergency Power Supply}
Several approaches leveraging the V2X technology for emergency power supply have been proposed. 
Xu et al. \cite{EPS1} proposed a (dis-) charge scheduling approach assuming that each household possesses an EV, which moves back and forth between the household and a charge station/another household. 
Yang et al. \cite{EPS2} consider public EV participants and have proposed a two-stage approach that determines the participants' discharge scheduling and rewards (i.e., money). In contrast to these works, EVRP-EPS is a routing problem (not discharge time scheduling) and ignores discharge rewarding as we here assume that maintainers of the base stations (our company's workers) operate EVs.

Although it does not use EVs, the approach in \cite{EVRP-EPS-close} shares our motivation that vehicles supply electricity to the affected building during/after a disaster. It uses power supply vehicles to supply electricity to batteries equipped with shelters.
The authors formulated this problem as a variant of VRPTW, where the goal is to minimize the sum of arrival time and waiting time for each shelter. On the other hand, EVRP-EPS is a variant of IRP, where there is no explicit time window, and multiple visits to the same node are allowed. This is because base station batteries are depleted relatively quickly and require multiple supplies during a time horizon.



\section{Conclusion and Future work}
In this paper, we formulated a base station relief by EVs as a variant of the EVRPs and proposed a solver that combines a rule-based vehicle selector and an RL-based node selector. We compared our solver with baselines on both synthetic datasets and real datasets. The results show that our solver consistently outperforms baselines in terms of the objective value and the computation time. We also analyzed the scalability and generalization performance, demonstrating our solver's capability for large-scale problems and unseen situations.

On the other hand, our solver currently faces two limitations: unbalanced work among EVs and only considering Euclidean distance.
Balancing the amount of work (e.g., the travel distance and the count of visiting nodes) among EVs is critical for improving their fairness if the drivers are humans.
Moreover, the effectiveness of our method would be reduced if the actual distance between two points differs significantly from the Euclidean distance.
In future work, we will address these problems by leveraging route-balancing techniques and a variant of Transformers that considers edge features, representing the actual travel distance.

\section*{Acknowledgement}
In the illustrations, we used free materials from Flaticon.com (\url{https://www.flaticon.com/}).


\bibliographystyle{ACM-Reference-Format} 
\bibliography{main}


\end{document}